\documentclass[draft,11pt,english]{article}

\usepackage{amsfonts,amsmath,amsthm,amscd,amssymb,latexsym,cite,verbatim,texdraw,floatflt,
caption2}

\theoremstyle{plain}
\newtheorem{theorem}{Theorem}[section]

\theoremstyle{question}

\newcommand{\keywords}{\textbf{Key words. }\medskip}
\newcommand{\subjclass}{\textbf{MSC }\medskip}
\renewcommand{\abstract}{\textbf{Abstract. }\medskip}

\numberwithin{equation}{section}

\begin{document}

\title{On hyperballeans of bounded geometry}

\author{Igor Protasov and Ksenia Protasova}



\maketitle

\begin{abstract}
A ballean (or coarse structure) is a set endowed with some family of subsets, the balls, is such a way that balleans with corresponding morphisms can be considered as asymptotic counterparts of uniform topological spaces. For a ballean $\mathcal{B}$ on a set $X$, the hyperballean $\mathcal{B}^{\flat}$ is a  ballean naturally defined on the set $X^{\flat}$ of all bounded subsets of $X$. We describe all balleans with hyperballeans of bounded geometry and analyze the structure of these hyperballeans.
\end{abstract}



\section{Introduction and preliminaries}

Following  [7], [8], we say that a {\it ball structure} is a triple $\mathcal{B}=(X, P, B)$, where $X, P$ are non-empty sets, and for all $x\in X$  and $\alpha\in P$, $B(x, \alpha)$ is a subset of $X$ which is called a {\it ball of radius} $\alpha$ around $x$. It is supposed that $x\in B(x, \alpha)$ for all $x\in X$, $\alpha\in P$. The set $X$ is called the {\it support} of $\mathcal{B}$ $P$ is called the {\it set of radii}.

Given any $x\in X$, $A\subseteq X$, $\alpha\in P$, we set $$B^{\ast}(x, \alpha)=\{y\in X:  x\in B(y,\alpha)\}, \  \  B(A, \alpha)=\bigcup_{a\in A}B(a,\alpha).$$

A ball structure $\mathcal{B}=(X, P, B)$  is called a {\it ballean} if

$\bullet$ for any $\alpha, \beta\in P$, there exist $\alpha^{\prime}, \beta^{\prime}\in P$ such that, for every $x\in X$, $$B(x,\alpha)\subseteq  B^{\ast}(x,\alpha^{\prime}),  \  \   B^{\ast}(x,\beta)\subseteq  B(x,\beta^{\prime});$$

$\bullet$ for every $\alpha, \beta\in P$, there exist $\gamma\in P$  such that, for every $x\in X$, $$B(B(x,\alpha), \beta)\subseteq  B(x,\gamma);$$

$\bullet$ for any $x,y\in X$, there exists $\alpha\in P$ such that $y\in B(x,\alpha)$.
\vskip 5pt

We note that a ballean can be considered as an asymptotic counterpart of a uniform space, and could be defined [9] in terms of entourages of the diagonal $\Delta_{X}$  in $X\times X$. In this case a ballean is called a {\it coarse structure}. For categorical look at the ballean and coarse structures as "two faces of the same coin" see [4].

Let $\mathcal{B}=(X, P, B)$,  $\mathcal{B^{\prime}}=(X^{\prime}, P^{\prime}, B^{\prime})$ be balleans. A mapping $f: X\longrightarrow X^{\prime}$ is called {\it coarse} if, for every $\alpha\in P$, there exists $\alpha\prime\in P^{\prime}$ such that
$$f(B(x,\alpha))\subseteq  B^{\prime}(f(x), \alpha^{\prime}).$$

A bijection $f: X\longrightarrow X^{\prime}$ is called an {\it asymorphism} between $\mathcal{B}$
and $\mathcal{B}^{\prime}$  if $f$ and $f^{-1}$ are coarse mappings.
In this case $\mathcal{B}$
and $\mathcal{B}^{\prime}$ are called {\it asymorphic}. If $X=X^{\prime}$ and the identity mapping $id: X\longrightarrow X^{\prime}$ is an asymorphism, we identify $\mathcal{B}$
and $\mathcal{B}^{\prime}$, and write $\mathcal{B}=\mathcal{B}^{\prime}.$
Given any ballean $\mathcal{B}=(X, P, B)$, replacing each ball $B(x,\alpha)$ to $B(x,\alpha)\cap B^{\ast}(x,\alpha)$, we  get the same ballean, so in what follows  we suppose that $B(x,\alpha)= B^{\ast}(x,\alpha)$.

Let $\mathcal{B}=(X, P, B)$ be a ballean. Each non-empty subset $Y$ of $X$ defines a {\it subballean} $\mathcal{B}_{Y}=(X, P, B_{Y})$, where $B_{Y}(y,\alpha)=Y\cap B(y,\alpha)$.
A subset $Y$ is called {\it large} if $X=B(Y,\alpha)$ for some $\alpha\in P$. Two  balleans $\mathcal{B}$ and $\mathcal{B}^{\prime}$ with the support $X$ and $X^{\prime}$ are called {\it coarsely equivalent} if there exist large subsets $Y\subseteq X$ and $Y^{\prime}\subseteq X^{\prime}$ such that the balleans $\mathcal{B}_{Y}$ and $\mathcal{B}^{\prime}_{Y^{\prime}}$ are asymorphic.

For a ballean  $\mathcal{B}=(X, P, B)$, a subset $Y$ of $X$ is called {\it bounded} if there exist $x\in X$ and $\alpha\in P$ such that $Y\subseteq B(x,\alpha)$.  A ballean $\mathcal{B}$ is called {\it bounded} if the support $X$ is bounded. Each bounded ballean is coarsely equivalent to a ballean whose support is a singletone.

Now we are ready to introduce the main subject of the note. For a ballean  $\mathcal{B}=(X, P, B)$, we denote by $X^{\flat}$ the family of all non-empty bounded subsets of, consider the ballean $\mathcal{B}^{\flat}=(X^{\flat}, P, B^{\flat})$, where
$$ B^{\flat}(Y,\alpha)=\{Z\in X^{\flat}: Z\subseteq B(Y,\alpha), Y\subseteq B(Z, \alpha)\},$$
and say that $\mathcal{B}^{\flat}$ is the {\it hyperballean} of $\mathcal{B}$.

For $\alpha\in P$, a subset $S$ of $X$ is called $\alpha$-discrete if $B(x,\alpha)\cap S=\{x\}$ for each $x\in S$. We say that $\mathcal{B}$ is of {\it bounded geometry} if there exist $\alpha\in P$ and a function $f: P\longrightarrow \mathbb{N}$ such that if $S$ is an $\alpha$-discrete subset of a ball $B(x,\beta)$ then $\mid \delta\mid\leq f(\beta)$. A ballean $\mathcal{B}$ is called {\it uniformly locally finite} if, for every  $\beta\in P$, there is $n(\beta)\in \mathbb{N}$ such that $\mid B(x,\beta)\mid\leq n(\beta)$ for every $x\in X$. By [6], $\mathcal{B}$ is of bounded geometry if and only if there exists large subset $Y$ of $X$ such that $\mathcal{B} _{Y}$  is uniformly locally finite.

It should be mentioned that the notion of bounded geometry went from asymptotic topology where metric spaces of bounded geometry play the central part [5]. For interrelations between balleans of bounded geometry and $G$-spaces see [6].

Every metric space $(X,d)$ defines the {\it metric ballean } $(X, \mathbb{R}^{+}, B_{\alpha})$, where $B_{d}(x,r)=\{y\in X: d(x,y)\leq r\}$. A ballean $\mathcal{B}$ is called {\it metrizable} if $\mathcal{B}$ is asymorphic to some metric ballean. By [ 8, Theorem 2.1.1], for a ballean  $\mathcal{B}$, the following statements are equivalent: $\mathcal{B}$ is metrizable, $\mathcal{B}$ is coarsely equivalent to some metrizable ballean, the set $P$ has a countable confinal subset $S$. We recall $S$ is confinal if, for every $\beta\in P$ there is $\alpha\in S$ such that $\alpha>\beta$. Here $\alpha>\beta$ means that $B(x,\beta)\subseteq B(x,\alpha)$ for each $x\in X.$ Applying this criterion, we conclude that, for every metrizable ballean $\mathcal{B}$, the hyperballean $\mathcal{B}^{\flat}$ is metrizable.

\section{Results }

For a non-empty set $X$ and the family $\mathcal{F}_{X}$ of all finite subsets of $X$, we denote by ${\bf F}_{X}$ the ballean $(X, \mathcal{F}_{X}, B_{{\bf F}})$ where

\[
B_{{\bf F}}(x,F)=
\begin{cases}
\{x\} & \text{if $x\notin F$;} \\
F & \text{if $x\in F$.}
\end{cases}
\]
Then ${\bf F}_{X}^{ \flat}  = \{\mathfrak{F}_{X}\backslash \{\emptyset\}, \mathfrak{F}_{X},
B_{{\bf F}} ^{ \flat}\}$, where $B_{{\bf F}} ^{ \flat}(H, F)=\{H\}$ if $H\cap F=\emptyset$ and
$B_{{\bf F}}(H, F)= \{(H\backslash F)\cup Z: Z\subseteq H, Z\neq \emptyset\}$ otherwise.

The ballean ${\bf F}_{\omega}, \omega=\{0,1\ldots\}$ is metrizable (say, by the metric $d(m,n)=|2^{m} -2^{n}|$), so ${\bf F}_{\omega}^{\flat}$ is also metrizable (say, by the Hausdorff metric $\flat_{H}$). At the end of the note, we point out some more explicit metrization of ${\bf F}_{\omega}^{\flat}$.

\begin{theorem}


For an unbounded ballean $\mathcal{B}=(X, P, B)$, the following statements hold:

$(i)$ $\mathcal{B}^{\flat}$ is uniformly locally finite if and only if $\mathcal{B}= {\bf F} _{X}$;

$(ii)$ $\mathcal{B}^{\flat}$ is of bounded geometry if and only if there exists a large subset $Y$ of $X$ such that $\mathcal{B}_{Y}= {\bf F} _{Y}$.

\end{theorem}

For a cardinal $\kappa$, we denote by ${\bf Q} _{\kappa}$ the ballean with the support
$${\bf Q} _{\kappa}=\{(x_{\alpha})_{\alpha<\kappa} : x_{\alpha}\in\{0,1\},  x_{\alpha}=0 $$
$$ \text{for   all   but   finitely  many} \ \ \alpha<\kappa \}, $$
the set of radii $\mathcal{F}_{\kappa}$ and the balls $$B_{{\bf Q}}((x_{\alpha})_{\alpha<\kappa}, F)= \{(y_{\alpha})_{\alpha<\kappa}: x_{\alpha} = y_{\alpha}  \ \ \text{  for  all} \ \ \alpha\in\kappa\setminus F \}.$$
The ballean ${\bf Q}_{\omega}$ is known as the Cantor macrocube and sometimes is denoted by $2^{<\omega}$ or $2^{<\mathbb{N}}$. For characterization of balleans coarsely equivalent to the Cantor macrocube see [3]. In [1], $2^{<\kappa}$ denotes the ballean of all $\{0,1\}$  $\kappa$-sequences $(x_{\alpha})_{\alpha<\kappa}$ such that $|\{\alpha<\kappa: x_{\alpha}=1\}|<\kappa .$

A ballean $\mathcal{B}=(X,P,B)$ is called {\it asymptotically scattered} if, for every unbounded subset $Y$ of $X$, there is $\alpha\in P$, such that, for every $\beta\in P$, there exists
$y\in Y$ such that $$(B(y,\beta)\setminus B(y,\alpha))\bigcap Y=\emptyset . $$
For asymptotically scattered subbaleans of group balleans see [2].

For a ballean $\mathcal{B}=(X,P,B)$, the subset $Y, Z$ of $X$ are called {\it close} if there exists $\alpha\in P$ such that $Y\subseteq B(Z,\alpha)$, $Z\subseteq B(Y,\alpha)$.


\begin{theorem}
Let $\kappa$ be an infinite cardinal,
$n\in \mathbb{N}$, $[\kappa]^{n}= \{F\subset \kappa: |F|=n\}$, $x\in\kappa$. Then the following statements hold:

$(i)$  the subbalean of  ${\bf F}_{\kappa}^{\flat}$ with the support $[\kappa]^{n}$ is asymptotically scattered;

 $(ii)$ the subbalean of  ${\bf F}_{\kappa}^{\flat}$ with the support $\{F\in \mathfrak{F}_{\kappa} :x\in F\}$ is asymorphic to ${\bf Q}_{\kappa}$;

 $(iii)$ ${\bf F}_{\omega}^{\flat}$ can be partitioned into countably many pairwise close Cantor macrocubes but ${\bf F}_{\omega}^{\flat}$ is not coarsely equivalent to ${\bf Q}_{\omega}$.
 \end{theorem}

 At the end of the note, we describe some explicit asymorphic embedding of ${\bf F}_{\omega}^{\flat}$  into ${\bf Q}_{\omega}$.

\section{Proofs}

{\it Proof of Theorem 1.2.} $(i)$ By the definition of balls in  ${\bf F}_{\omega}^{\flat}$,
${\bf F}_{\omega}^{\flat}$ is uniformly locally finite.

If the identity mapping $id: X\longrightarrow X$  is not an asymorphism between $\mathcal{B}$ and ${\bf F}_{X}$ then we can choose $\alpha\in P$ and a sequence $(x_{n})_{n<\omega}$ in $X$ such that $|B(x_{n}, \alpha)|>1$ and $B(x_{i},\alpha)\bigcap B(x_{j},\alpha)=\emptyset$ for all $i<j<\omega$. For each $i<\omega$, we pick $y_{i}\in B(x_{i}, \alpha)$, $y_{i}\neq x_{i}$, put $X_{n}=\{x_{0}, \ldots , x_{n}\}$, $X_{n,i}=X_{n}\bigcup\{y_{i}\}$, $i\leq n<\omega$. Then $X_{n,i}\in B^{\flat}(X_{n}, \alpha)$, so $|B^{\flat}(X_{n}, \alpha)|>n$ and $\mathcal{B}$ is not uniformly locally finite.
\vskip 5pt

$(ii)$  We assume that $Y$ is a large subset of $X$  and choose $\beta\in P$ such that $B(Y,\beta)=X$. For each $x\in X$, we pick $y_{x}\in Y$ such that $y_{x}\in B(x,\beta)$. If $F\in X^{\flat}$ then $\{y_{x}: x\in F\}\in Y^{\flat}$ and $F\in B^{\flat}(\{y_{\kappa}: x\in F\},\beta)$. It follows that
$B^{\flat}(Y^{\flat}, \beta)=X^{\flat}$, $Y^{\flat}$  is large in $X^{\flat}$ so $\mathcal{B}^{\flat}_{Y}$ and $\mathcal{B}^{\flat}$ are coarsely equivalent. In particular, if $\mathcal{B}_{Y}={\bf F}_{Y}$, we conclude that $\mathcal{B}$ is of bounded geometry.

We suppose that  $\mathcal{B}^{\flat}$ is of bounded geometry and let $\alpha\in P$ and $f: P\longrightarrow \mathbb{N}$ witness this property. Using Zorn's lemma, we choose a maximal by inclusion subset $Y$ of $X$ such that $B(y,\alpha)\bigcap B(y^{\prime},\alpha)=\emptyset$ for all distinct $y, y^{\prime}\in Y$. We show that $\mathcal{B}_{Y}={\bf F}_{Y}$.

If the identity mapping $id: Y\longrightarrow Y$ is not an asymorphism between $\mathcal{B}_{Y}$ and ${\bf F}_{Y}$ then there are $\beta\in P$  and a sequence $(y_{n})_{n\in\omega}$  in $Y$ such that $|B_{Y}(y_{n}, \beta)|> 1$ and $B_{Y}(y_{i}, \beta) \bigcap B_{Y}(y_{j}, \beta)=\emptyset$ for all $i<j<\omega$. For each $i<\omega$, we pick $z_{i}\in B_{Y}(y_{i}, \beta)$, $z_{i}=y_{i}$, put
$Y_{n}=\{y_{0}, \ldots , y_{n}\} $, $Y_{n,i} = Y_{n} \bigcup \{y_{i}\}$, $i\leq n<\omega$. Then $Y_{n,i}\in B^{\flat}(Y_{n}, \beta)$ and the set $\{Y_{n,i}: i\leq n\}$ is $\alpha$-discrete. Thus, for $n> f(\beta)$ we get a contradiction with the choice of $\alpha$ and  $f$.  $ \ \ \ \ \Box$

{\it Proof of Theorem 2.2.} $(i)$  We say that a subset of a ballean is {\it asymptotically scattered} if corresponding subballean has this property. We use the following observation: the union of two asymptotically scattered subsets is asymptotically scattered (see [2]).

We note that every unbounded subset in ${\bf F}_{\kappa}^{\flat}$ is infinite and proceed on induction by $n$. For $n=1$, the statement is evident: given any $H\in \mathfrak{F}_{\kappa}$ and an infinite subset $Y$ of $[\kappa]^{1}$, we take $\{y\}\in Y$, $y\in H$ and get $B_{{\bf F}} (\{y\}, H)=\{y\}$.

Assuming that the statement is true for $[\kappa]^{n}$, let $Y$ be an infinite subset of $[\kappa]^{n+1}$.  For each $F\in [\kappa]^{n+1}$, we denote by $\min F$ and $\max F$, the minimal and maximal elements of $F$ with respect to the ordinal ordering of $\kappa$  and consider two cases.

Case: the set $\{\min F: F\in Y\}$ is infinite. We take an arbitrary $H\in \mathfrak{F}_{\kappa}$ and choose $F\in Y$ such that $\max H< \min F$. Then $B_{{\bf F}}^{\flat}(F,H)=\{F\}$.

Case: the set $\{\min F: F\in Y\}$ is finite, $\{\min F: F\in Y\}= x_{1}, \ldots , x_{n}$. For each $i\in \{1,\ldots, n\}$, we denote $Z_{i}=\{F\in[\kappa]^{n+1}: x_{i}\in F\}$. We note that $Z_{i}$
is asymorphic to $[\kappa]^{n}$ and, by the inductive assumption, $Z_{i}$ is asymptotically scattered.

Then $Z_{1}\bigcup \ldots \bigcup Z_{n}$ is asymptotically scattered, $Y\subseteq Z_{1}\bigcup \ldots \bigcup Z_{n}$  and we can use definition of asymptotically scattered subsets to choose $\alpha\in \mathfrak{F}_{X}$ suitable for $Y$.

$(ii)$  We use the standard bijection $\chi: \mathfrak{F}_{\kappa}\longrightarrow Q_{\kappa}$
 defined by $\chi(K)= (x_{\alpha})_{\alpha<\kappa}$, where $x_{\alpha}=1$ if and only  if $\alpha\in K$. Then the restriction of $\chi$ to $\{F\in \mathfrak{F}_{\kappa}\}: x\in F$ is a asymorphic  embedding. Indeed, to verity this property we may use as radii in ${\bf F}^{\flat}$ only balls containing $x$. Clearly, $\chi\{F\in \mathfrak{F}_{\kappa}: \kappa\in F\}$ is asymorphic to ${\bf Q}_{\kappa}$.   $ \ \ \ \Box$

 $(iii)$  For every $n\in\omega$, let $\mathcal{M}_{n}=\{F\in \mathfrak{F}_{X}: \min F=n\}$.
 Applying $(ii)$, we see that $\mathcal{M}_{n}$ is
asymorphic to ${\bf Q}_{\omega}$. We take arbitrary $i,j\in \omega$, denote $m=\max \{i,j\}$,
$I_{m}=\{0,\ldots,m\}$. Then $\mathcal{M}_{i}\subseteq B_{{\bf F}}^{\flat}(\mathcal{M}_{j}, I_{m})$, $\mathcal{M}_{j}\subseteq B_{{\bf F}}^{\flat}(\mathcal{M}_{i}, I_{m})$ so
$\mathcal{M}_{i}$, $\mathcal{M}_{j}$ are close.

Given $H\in \mathfrak{F}_{\omega}$, we take $F\in \mathfrak{F}_{\omega}$ such that $\max H< \min F$. Then $B_{{\bf F}}^{\flat}(F,H)=\{F\}$.
In terminology of [3], it means that ${\bf F}^{\flat}_{\omega}$ has an asymptotically isolated balls but
every ballean coarsely  equivalent to ${\bf Q}_{\omega}$ has no isolated balls.   $ \  \ \  \Box$

To embed asymorphically  ${\bf F}^{\flat}_{\omega}$ into ${\bf Q}_{\omega}$, we use $2\mathbb{N}$ in place of $\omega$. We define a mapping $f: \mathfrak{F} _{2\mathbb{N}}\ \{\emptyset\}\longrightarrow Q_{\omega}$ by the $f(K)=(x_{n})_{n<\omega}$, where $x_{n}=1$ if and only if $n\in\{\min K - 1\}\bigcup K$. We note that the set $S=f(\mathfrak{F} _{2\mathbb{N}} \ \setminus \{\emptyset\})$ consists of all sequences $(x_{n})_{n<\omega}$ with at  least two
non-zero coordinates and such that the first non-zero coordinate of $(x_{n})_{n<\omega}$ is odd and all other are even.  For each $K\in \mathfrak{F}_{2\mathbb{N}}\backslash \{\emptyset\}$ and $n\in \mathbb{N}$, we have $$f(B_{{\bf F}}(K, \{2,4,\ldots,2n\})= S\bigcap B_{{\bf Q}}(f(K), \{1,2.\ldots,2n\}),$$
witnessing that $f$ is an asymorphic embedding of ${\bf F}^{\flat}_{2\mathbb{N}}$ into  ${\bf Q}_{\omega}$.

With this representation, ${\bf F}^{\flat}_{\omega}$ can be easily metrizable by means of restriction to $S$ of the stadard metric $d$ on $Q_{\omega}: d((x_{n})_{n\in\omega}, (y_{n})_{n\in\omega})=\min\{m: x_{n}=y_{n} \ \ \text{ for all} \ \ n\geq m\} $.

\
\begin{description}

\item{1.} T. Banakh, I. Protasov, D. Repov$\breve{s}$, S. Slobodianiuk, \emph{Classifying homogeneous celular ordinal balleans up to coarse equivalence},
preprint (arxiv: 1409.3910v2).

 \item{2.} T. Banakh, I. Protasov,  S. Slobodianiuk, \emph{Scattered subsets of groups}, Ukr. Math. Zh. \textbf{67} (2015), P. 304-312.

\item{3.} T. Banakh, I. Zarichnyi, \emph{Characterizing the Cantor bi-cube in asymptotic categories}, Groups, Geometry and Dynamics \textbf{5} (2011), P. 691-728.

\item{4.} D. Dikranjan, N. Zava,\emph{ Some categorical aspects of coarse spaces and balleans,} Topology Appl., to appear.

\item{5.} A. Dranishnikov, \emph{Asymptotic Topology,} Russian Math. Surveys, \textbf{55} (2000), 1085-1129.

 \item{6.} I. Protasov,\emph{ Balleans of bounded geometry and $G$-spaces,}  Math. Stud.  \textbf{30} (2008), 61-66.

\item{7.} I. Protasov, T. Banakh, \emph{ Ball Structures and
Colorings of Graphs and Groups}, Math. Stud. Monogr. Ser. 11, Vol. 11,
VNTL Publisher, Lviv, 2003.

\item{8.} I. Protasov, I. Zarichnyi, \emph{General Asymptology}, Math. Stud. Monogr. Ser 12, Vol. 11, VNTL Publisher, Lviv, 2007.

\item{9.} J. Roe,\emph{Lectures on coarse geometry,} Univ. Lecture Series 31, American Mathematical Society, Providence, 2003.

\end{description}
\bigskip

\bigskip

CONTACT INFORMATION

\medskip

Igor Protasov \\
Department of Cybernetics  \\
         Kyiv University\\
         Prospect  Glushkova 2, corp. 6  \\
         03680 Kyiv, Ukraine \\ e-mail: i.v.protasov@gmail.com

\medskip

Ksenia Protasova\\
Department of Cybernetics  \\
        Kyiv University\\
         Prospect  Glushkova 2, corp. 6  \\
         03680 Kyiv, Ukraine \\ e-mail: k.d.ushakova@gmail.com

\end{document}